\documentclass[12pt,oneside,draft]{amsart}
\usepackage{amssymb,amsmath}
\usepackage{graphpap}

 
\newtheorem{theorem}[subsection]{Theorem}
\newtheorem*{theorem*}{Theorem}
\newtheorem{lemma}[subsection]{Lemma}

\newtheorem{corollary}[subsection]{Corollary}

\theoremstyle{definition}
\newtheorem{definition}[subsection]{Definition}

\theoremstyle{remark}
\newtheorem{remark}[subsection]{Remark}

\newcommand{\mt}[1]{\operatorname{#1}}
\newcommand{\EEE}{{\mathbb E}}
\newcommand{\DDD}{{\mathbb D}}
\newcommand{\AAA}{{\mathbb A}}
\newcommand{\QQ}{{\mathbb Q}}
\newcommand{\ZZ}{{\mathbb Z}}
\newcommand{\CC}{{\mathbb C}}
\newcommand{\OO}{{\mathcal O}}
\newcommand{\RR}{{\mathbb R}}
\newcommand{\PP}{{\mathbb P}}
\newcommand{\FF}{{\mathcal F}}
\newcommand{\TT}{{\mathcal T}}
\newcommand{\NN}{{\mathbb N}}
\newcommand{\FFF}{{\mathbb F}}
\newcommand{\Supp}{\mt{Supp}}
\newcommand{\Sing}{\mt{Sing}}
\newcommand{\Diff}{\mt{Diff}}
\newcommand{\rank}{\mt{rank}}
\newcommand{\ord}{\mt{ord}}
\newcommand{\wt}{\mt{wt}}
\newcommand{\D}{{\Delta}}
\newcommand{\Hom}{\mt{Hom}}
\newcommand{\Exc}{\mt{Exc}}
\newcommand{\codim}{\mt{codim}}
\newcommand{\mult}{\mt{mult}}
\newcommand{\Spec}{\mt{Spec}}
\newcommand{\Gal}{\mt{Gal}}

\newcommand{\down}[1]{\llcorner #1 \lrcorner}
\newcommand{\up}[1]{\ulcorner #1 \urcorner}
\newcommand{\fr}[1]{\{ #1\}}

\title{Classification of log Enriques surfaces with {\large $\delta=1$}}
\author{S.~A.~Kudryavtsev}

\date{}

\address{Department of Algebra, Faculty of Mathematics,
Moscow State Lomonosov University, 117234 Moscow, Russia}

\email{kudryav@mech.math.msu.su}

\begin{document}
\begin{abstract} Log Enriques surfaces with $\delta=1$ are classified.
\end{abstract}
\maketitle

\section*{\bf {Introduction}}
Let $S$ be a projective surface with klt singularities and with numerically trivial
canonical class $K_S\equiv 0$. Then $S$ is called {\it a log Enriques surface}.
\par
Let us consider the following invariant
\begin{eqnarray*}
\delta(S)=\#\Big\{E| E \ \text{is exceptional divisor with discrepancy}\\
a(E,0)\le -\frac67   \Big\}.
\end{eqnarray*}

By theorem \cite[5.1]{Sh2} we have $0\le \delta(S) \le 2$.
In this paper the classification of such surfaces with $\delta=1$ is given.
The classification of log Enriques surfaces with
$\delta=2$ was obtained in \cite{Kud4}.
\par
The classification method for surfaces with
$\delta=1$ doesn't distinguishes from one in the case
$\delta=2$. Moreover, the classification of exceptional log
Del Pezzo surfaces easily implies the classification of
log Enriques surfaces.
\par
The log Enriques surfaces often appear in many problems.
Such of them are the following ones:
the study of
Calabi--Yao varieties, the study of $K3$ surfaces
(in particular, see section \ref{K3}) and the problem of inductive
classification of strictly log canonical singularities.
\par
In the latter problem log Enriques surfaces can be realized as the exceptional
divisors of purely log terminal blow-ups of three-dimensional strictly
log canonical singularities.
\par
The research was partially supported by a grant 02-01-00441 from
the Russian Foundation of Basic Research and a grant INTAS-OPEN
2000\#269.

\section{\bf {Preliminary facts and results}}

All varieties are algebraic and are assumed to be defined over
$\CC$, the complex number field.
The main definitions, terminology and notations used in the paper are
given in \cite{Koetal}, \cite{PrLect}.

\par
The general hypersurface of degree $d$ in $\PP$ is denoted by $X_d$.
In many cases it is enough to require the irreducibility and reducibility of $X_d$.
In the case $p_a(C)=1$ the curve $X_d$ must have an ordinary double point.
More exactly about a structure of
$X_d$ see remark 2.2 \cite{Kud3}.

\begin{theorem}\label{ldmain}
Let $S$ be a log Enriques surface with $\delta=1$. Denote the exceptional
curve with discrepancy $a(\widetilde{C},0)=-a\le -\frac67$ by
$\widetilde{C}$. Let $f\colon \widetilde S\to  S$ be an extraction of
$\widetilde{C}$
$($i.e. $f$ is a birational contraction such that $\Exc f=\widetilde C)$.
\\
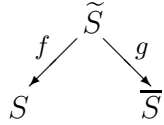
\begin{center}
\begin{picture}(80,46)(0,0)
\put(14,10){$S$}
\put(40,40){\vector(-1,-1){18}}
\put(42,42){$\widetilde S$}
\put(50,40){\vector(1,-1){18}}
\put(64,10){$\overline{S}$}
\put(24,34){\footnotesize{$f$}}
\put(62,34){\footnotesize{$g$}}
\end{picture}
\end{center}
\par
Then there exists a birational contraction
$g\colon \widetilde S\to  \overline{S}$ with the following properties:
$\rho(\overline{S})=1$, $g$
doesn't contract the curve $\widetilde{C}$. Put $C=g(\widetilde{C})$.
Then the pair $\big(\overline{S},aC\big)$
is one of the following ones:
\begin{enumerate}
\item  $\overline{S}=\PP(1,2,3)$ and $aC=\frac67X_7$.
It is the case $6-2(\rm ell)$.
\item  $\overline{S}=\PP(1,3,4)$ and $aC=\frac89X_9$.
It is the case $8-1(\rm ell)$.
\item  $\overline{S}=\PP(1,3,5)$ and $aC=\frac9{10}X_{10}$.
It is the case $9-1(\rm ell)$.
\item  $\overline{S}=\PP(1,5,8)$ and $aC=\frac78X_{16}$.
It is the case $22-1(\rm ell)$.
\item The minimal resolution of $\overline{S}$ is one of the following ones.

\begin{center}
\begin{picture}(240,70)(0,0)
\put(18,16){\circle*{8}}
\put(0,13){\tiny{$-2$}}
\put(18,20){\line(0,1){8}}
\put(18,32){\circle*{8}}
\put(0,29){\tiny{$-2$}}
\put(18,36){\line(0,1){8}}
\put(18,48){\circle*{8}}
\put(0,45){\tiny{$-3$}}
\put(18,52){\line(0,1){8}}
\put(18,64){\circle*{8}}
\put(0,61){\tiny{$-2$}}
\put(22,48){\line(1,0){8}}
\put(38,48){\circle{16}}
\put(32,58){\tiny{$-1$}}
\put(46,48){\line(1,0){8}}
\put(58,48){\circle*{8}}
\put(52,54){\tiny{$-2$}}
\put(62,48){\line(1,0){8}}
\put(70,44){\fbox{$C$}}
\put(72,60){\tiny{$+5$}}
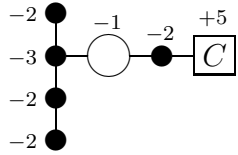
\put(15,0){\tiny{Case $51-2({\rm ell)}$}}
\put(162,16){\circle*{8}}
\put(144,13){\tiny{$-2$}}
\put(162,20){\line(0,1){8}}
\put(162,32){\circle*{8}}
\put(144,29){\tiny{$-2$}}
\put(162,36){\line(0,1){8}}
\put(162,48){\circle*{8}}
\put(144,45){\tiny{$-3$}}
\put(162,52){\line(0,1){8}}
\put(162,64){\circle*{8}}
\put(144,61){\tiny{$-2$}}
\put(166,48){\line(1,0){8}}
\put(178,48){\circle*{8}}
\put(172,54){\tiny{$-2$}}
\put(182,48){\line(1,0){8}}
\put(194,48){\circle{8}}
\put(188,54){\tiny{$-1$}}
\put(198,48){\line(1,0){8}}
\put(210,48){\circle*{8}}
\put(204,54){\tiny{$-3$}}
\put(214,48){\line(1,0){8}}
\put(222,44){\fbox{$C$}}
\put(224,60){\tiny{$+5$}}
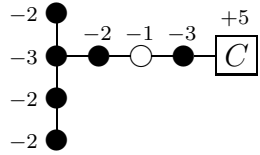
\put(160,0){\tiny{Case $51-6({\rm ell)}$}}
\end{picture}
\end{center}
We have $a=\frac{10}{11}$ and $a=\frac78$ in the cases
$51-2({\rm ell)}$ and $51-6({\rm ell)}$ respectively.

\item  $\overline{S}=\PP(2,3,7)$ and $aC=\frac67X_{14}$.
It is the case $52-2(\rm ell)$.

\item The minimal resolution of $\overline{S}$ is one of the following ones.
\begin{center}
\begin{picture}(200,70)(0,0)
\put(8,16){\circle*{8}}
\put(-10,13){\tiny{$-2$}}
\put(8,20){\line(0,1){8}}
\put(8,32){\circle*{8}}
\put(-10,29){\tiny{$-2$}}
\put(8,36){\line(0,1){8}}
\put(8,48){\circle*{8}}
\put(-10,45){\tiny{$-2$}}
\put(8,52){\line(0,1){8}}
\put(8,64){\circle*{8}}
\put(-10,61){\tiny{$-2$}}
\put(12,48){\line(1,0){8}}
\put(24,48){\circle*{8}}
\put(18,54){\scriptsize{$-3$}}
\put(28,48){\line(1,0){8}}
\put(40,48){\circle{8}}
\put(34,54){\tiny{$-1$}}
\put(44,48){\line(1,0){8}}
\put(56,48){\circle*{8}}
\put(50,54){\tiny{$-2$}}
\put(60,48){\line(1,0){8}}
\put(68,44){\fbox{$C$}}
\put(70,60){\scriptsize{$+4$}}
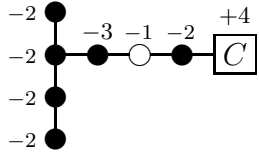
\put(5,0){\tiny{Case $53-2({\rm ell)}$}}
\put(140,22){\circle*{8}}
\put(134,10){\scriptsize{$-2$}}
\put(144,22){\line(1,0){16}}
\put(164,22){\circle*{8}}
\put(158,10){\scriptsize{$-3$}}
\put(168,22){\line(1,0){16}}
\put(188,22){\circle*{8}}
\put(182,10){\scriptsize{$-2$}}
\put(192,22){\line(1,0){16}}
\put(212,22){\circle*{8}}
\put(206,10){\scriptsize{$-2$}}
\put(216,22){\line(1,0){16}}
\put(236,22){\circle*{8}}
\put(230,10){\scriptsize{$-2$}}
\put(164,24){\line(0,1){8}}
\put(164,36){\circle{8}}
\put(146,33){\scriptsize{$-1$}}
\put(164,40){\line(0,1){8}}
\put(164,52){\circle*{8}}
\put(146,49){\scriptsize{$-2$}}
\put(168,52){\line(1,0){8}}
\put(176,48){\fbox{$C$}}
\put(179,63){\tiny{$+3$}}
\put(192,52){\line(1,0){8}}
\put(204,52){\circle{8}}
\put(199,59){\scriptsize{$-1$}}
\put(204,48){\line(3,-2){34}}
\put(208,52){\line(1,0){16}}
\put(228,52){\circle*{8}}
\put(223,59){\scriptsize{$-2$}}
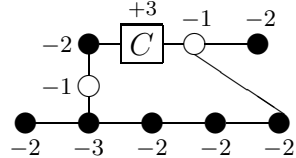
\put(175,0){\tiny{Case $54({\rm ell)}$}}
\end{picture}
\end{center}
We have $a=\frac89$ and $a=\frac67$ in the cases
$53-2({\rm ell)}$ and $54({\rm ell)}$ respectively.

\item  $\overline{S}=\PP(3,4,5)$ and $aC=\frac{12}{13}X_{13}$.
It is the case $18-1(+1)$.

\item  $\overline{S}=\PP(3,5,7)$ and $aC=\frac{15}{17}X_{17}$.
It is the case $25-1(+1)$.

\item The minimal of $\overline{S}$ is one of the following ones.

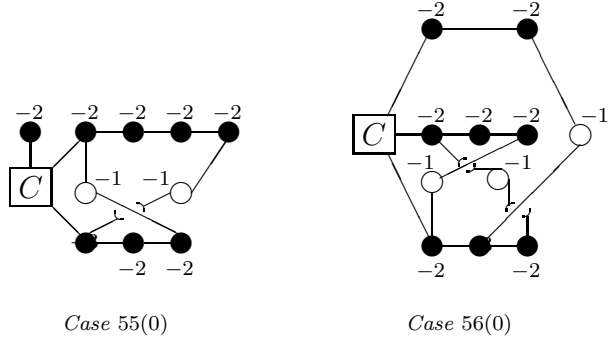
\begin{center}
\begin{picture}(285,125)(0,0)
\put(10,50){\fbox{$C$}}
\put(18,61){\line(0,1){10}}
\put(18,75){\circle*{8}}
\put(12,80){\tiny{$-2$}}
\put(26,61){\line(1,1){11}}
\put(39,75){\circle*{8}}
\put(33,80){\tiny{$-2$}}
\put(43,75){\line(1,0){10}}
\put(57,75){\circle*{8}}
\put(51,80){\tiny{$-2$}}
\put(61,75){\line(1,0){10}}
\put(75,75){\circle*{8}}
\put(69,80){\tiny{$-2$}}
\put(79,75){\line(1,0){10}}
\put(93,75){\circle*{8}}
\put(87,80){\tiny{$-2$}}
\put(26,47){\line(1,-1){11}}
\put(39,33){\circle*{8}}
\put(33,31){\tiny{$-3$}}
\put(43,33){\line(1,0){10}}
\put(57,33){\circle*{8}}
\put(51,21){\tiny{$-2$}}
\put(61,33){\line(1,0){10}}
\put(75,33){\circle*{8}}
\put(69,21){\tiny{$-2$}}
\put(39,52){\circle{8}}
\put(39,71){\line(0,-1){15}}
\put(42,55){\tiny{$-1$}}
\put(43,52){\line(2,-1){31}}
\put(75,52){\circle{8}}
\put(60,55){\tiny{$-1$}}
\put(71,52){\line(-2,-1){10}}
\put(58,45){\oval(6,6)[tr]}
\put(51,42){\line(-2,-1){12}}
\put(53,45){\oval(6,6)[bl]}
\put(79,52){\line(2,3){14}}
\put(30,0){\tiny{Case $55(0)$}}
\put(140,70){\fbox{$C$}}
\put(156,74){\line(1,0){10}}
\put(170,74){\circle*{8}}
\put(164,79){\tiny{$-2$}}
\put(174,74){\line(1,0){10}}
\put(188,74){\circle*{8}}
\put(182,79){\tiny{$-2$}}
\put(192,74){\line(1,0){10}}
\put(206,74){\circle*{8}}
\put(200,79){\tiny{$-2$}}
\put(170,32){\circle*{8}}
\put(174,32){\line(1,0){10}}
\put(164,20){\tiny{$-2$}}
\put(188,32){\circle*{8}}
\put(182,30){\tiny{$-3$}}
\put(192,32){\line(1,0){10}}
\put(206,32){\circle*{8}}
\put(200,20){\tiny{$-2$}}
\put(170,33){\line(-1,2){17}}
\put(170,114){\circle*{8}}
\put(164,119){\tiny{$-2$}}
\put(174,114){\line(1,0){28}}
\put(206,114){\circle*{8}}
\put(200,119){\tiny{$-2$}}
\put(170,115){\line(-1,-2){17}}
\put(226,74){\circle{8}}
\put(226,79){\tiny{$-1$}}
\put(206,114){\line(1,-2){18}}
\put(188,32){\line(1,1){38}}
\put(170,32){\line(0,1){20}}
\put(170,56){\circle{8}}
\put(160,61){\tiny{$-1$}}
\put(206,74){\line(-2,-1){33}}
\put(170,74){\line(1,-1){10}}
\put(183,62){\oval(6,6)[tl]}
\put(195,57){\circle{8}}
\put(197,61){\tiny{$-1$}}
\put(183,63){\oval(6,6)[br]}
\put(186,61){\line(1,0){10}}
\put(206,32){\line(0,1){11}}
\put(204,46){\oval(6,6)[br]}
\put(201,45){\oval(6,6)[tl]}
\put(199,57){\line(0,-1){10}}
\put(160,0){\tiny{Case $56(0)$}}
\end{picture}
\end{center}
We have $a=\frac{10}{11}$ and $a=\frac67$
in the cases $55(0)$ and $56(0)$ respectively.
\end{enumerate}
If $p_a(C)=1$ then the curve $C$ must have an ordinary double point.
\begin{proof}
This theorem immediately follows from theorem
2.1 \cite{Kud3}.
If we write ell in the brackets then
$p_a(C)=1$. If we write $q$ in the brackets then $p_a(C)=0$ and the self-intersection
index of proper transform of
$C$ on a minimal resolution of
$\overline{S}$ is equal to $q$.
The case $p_a(C)=0$ is rewritten from theorem
2.1 \cite{Kud3} completely. In the case $p_a(C)=1$ the curve $C$ must be rational.
Hence it must have an ordinary double point.
By lemma \ref{l1} the log surfaces with $p_a(C)=1$ are rewritten from theorem
2.1 \cite{Kud3} only if $a(\widetilde{C},0)=-\frac{k-1}{k}$, where
$k\in \ZZ_{\ge 2}$.
\end{proof}
\end{theorem}

\begin{lemma}\label{l1}
Let $(\CC^2,a\{xy=0\})$ be a klt pair. Then there exists a blow-up with
the discrepancy of exceptional divisor being equal to
0 if and only if
$a=\frac{k-1}{k}$, where $k\in \ZZ_{\ge 2}$.
\begin{proof} It is easy to prove that the required blow-up is toric,
i.e. it is a weighted blow-up of
$\CC^2$ with weights $(\alpha,\beta)$. Then
$a(E,0)=\alpha+\beta-1-(\alpha+\beta)a$, where $E$ is an exceptional divisor.
Hence $a(E,0)=0$ $\Longleftrightarrow$ $a=1-\frac1{\alpha+\beta}$.
\end{proof}
\end{lemma}

\begin{definition} Let $S$ be a log Enriques surface.
{\it A canonical index} $I=I(S)$ is called the following number
$I=\min\{n\in\ZZ_{>0}\mid nK_S\sim 0\}$. It is known that
$I\le 21$ \cite{Bl}, \cite{Z}.
\end{definition}

\begin{corollary} \label{index}
Let $S$ be a log Enriques surface with $\delta=1$. Then
$I\in\{7,8,9,10,11,13,17\}$.
\end{corollary}

So, the problem of classification of log Enriques surfaces with $\delta=1$
is to describe the following procedures.
At first we consider the extraction
$\widetilde S \to \overline{S}$ such that every exceptional divisor
$E_i$ has the discrepancy $a(E_i,aC)=0$.
Then we contract the proper transform of
$C$.
The number of such procedures is finite by the following easy fact.
\par
Let $(X,D)$ be a klt pair. Then the number of divisors
$E$ of the function field $\mathcal{K}(X)$ with
$a(E,D)\le 0$ is finite \cite[lemma 3.1.9]{PrLect}.

\subsection{} {\it Log Enriques surfaces and $K3$ surfaces.}\label{K3}
\par
Let $S$ be a log Enriques surface with $\delta=1$.
Consider its canonical cover
$\varphi\colon \widehat{S}=\Spec_{\OO_S}(\oplus^{I-1}_{i=0}\OO_S(-iK_S))\to S$, where
$I$ is an index of $S$. Since $I\ge 7$ (see corollary \ref{index}) then
$\widehat{S}$ is a $K3$ surface with at worst Du Val singularities
\cite{Bl}, \cite{Z} and
\begin{enumerate}
\item $\varphi$ is cyclic Galous cover of degree $I$, which is
etale over $S\backslash\Sing S$.
\item There exists a generator $g$ of $\Gal(\widehat{S}/S)\cong \ZZ_I$
such that $g^*\omega_{\widehat{S}}=\varepsilon_I\omega_{\widehat{S}}$, where
$\varepsilon_I=\exp(2\pi \sqrt{-1}/I)$ is a primitive root and
$\omega_{\widehat{S}}$ is a
nowhere vanishing regular
2-form on $\widehat{S}$.
\end{enumerate}
Let $\Delta(S)$ be an exceptional set of $\chi$, where $\chi$ is a minimal resolution
of
$\widehat{S}$. Then $\Delta(S)$ is a disconnected sum of divisors of Dynkin's type
$\AAA_i$, $\DDD_j$, $\EEE_k$.
So $\Delta(S)=(\oplus \AAA_{\alpha})\oplus(\oplus \DDD_{\beta})
\oplus(\oplus \EEE_{\gamma})$. Let us define
$\rank\Delta(S)=\sum\alpha+\sum\beta+\sum\gamma$.

\section{\bf {Classification of log Enriques surfaces with $\delta=1$ }}

Let $S$ be a log Enriques surface with $\delta=1$.
Then
$S$ can be constructed by the following way:
at first we extract some set $\TT$ of exceptional curves with discrepancy 0 for
the corresponding model of theorem \ref{ldmain}.
After
it we contract the proper transform
of $C$. So, the classification of log Enriques surfaces with $\delta=1$
is reduced to the description of sets $\TT$.
The curves from $\TT$ are enumerated by the natural numbers in the
corresponding figures.
\begin{definition}\label{def}
Let $T_i=|\TT_i|$ be a number of elements of set $\TT_i$,
$\underline{T_i}=\min\{t\mid t\in\TT_i\}$ and $\overline{T_i}=\max\{t\mid t\in\TT_i\}$.
\end{definition}

\begin{definition} Let $p_a(C)=1$. These are the cases (1)--(7)
from theorem \ref{ldmain}.
Then $\TT=\TT_1\cup\TT_2\cup\TT_3$, where
$\TT_1$ consists of the curves, which are contracted to the singular points
of $\overline{S}$ lying on
$C$,
$\TT_2$ consists of the curves, which are contracted to the singular point
of $C$, and
$\TT_3$ consists of the curves, which are contracted to the singular points
of $\overline{S}$ not lying on
$C$.
The set $\TT_2$ is considered up to symmetry.
For example, if $4\in \TT_2$ then $9\in \TT_2$ in the notation of theorem
\ref{first}. Note that $\TT_2\ne \emptyset$.
\end{definition}

\begin{theorem}\label{first}
In the case $6-2({\rm ell})$ the set $\TT$ must be one of the following ones.
\begin{enumerate}
\item Let $T_1=0$. Then $T_2\ge 2$ and $4\in \TT_2$.
\item Let $\TT_1=\{1\}$. Then $\{4,5\}\cap \TT_2\ne \emptyset$.
\item Let $\TT_1=\{3\}$. Then $4\in\TT_2$.
\item Let $\TT_1=\{2\}$ or $T_1\ge 2$. Then $\TT_2\ne \emptyset$.
\end{enumerate}
\end{theorem}

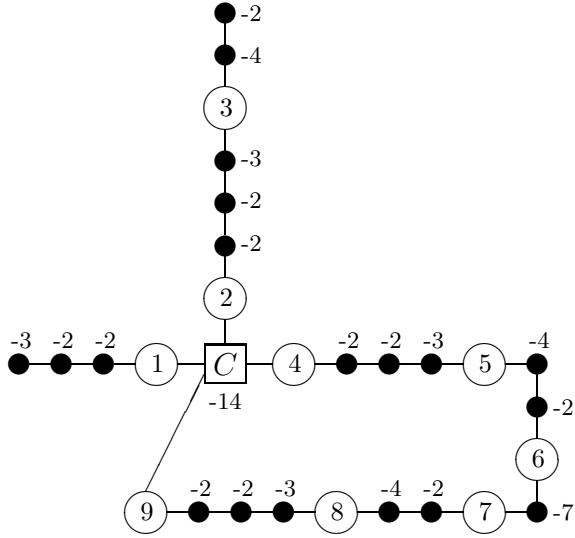
\begin{center}
\begin{picture}(210,190)(0,0)
\put(4,50){\circle*{8}}
\put(1,56){\scriptsize{-3}}
\put(8,50){\line(1,0){8}}
\put(20,50){\circle*{8}}
\put(17,56){\scriptsize{-2}}
\put(24,50){\line(1,0){8}}
\put(36,50){\circle*{8}}
\put(33,56){\scriptsize{-2}}
\put(40,50){\line(1,0){8}}
\put(56,50){\circle{16}}
\put(54,47){\footnotesize{1}}
\put(64,50){\line(1,0){10}}
\put(74,46){\fbox{$C$}}
\put(76,33){\scriptsize{-14}}
\put(82,57){\line(0,1){10}}
\put(82,75){\circle{16}}
\put(80,72){\footnotesize{2}}
\put(82,83){\line(0,1){8}}
\put(82,95){\circle*{8}}
\put(88,93){\scriptsize{-2}}
\put(82,99){\line(0,1){8}}
\put(82,111){\circle*{8}}
\put(88,109){\scriptsize{-2}}
\put(82,115){\line(0,1){8}}
\put(82,127){\circle*{8}}
\put(88,125){\scriptsize{-3}}
\put(82,131){\line(0,1){8}}
\put(64,50){\line(1,0){8}}
\put(82,147){\circle{16}}
\put(80,144){\footnotesize{3}}
\put(82,155){\line(0,1){8}}
\put(82,167){\circle*{8}}
\put(88,164){\scriptsize{-4}}
\put(82,171){\line(0,1){8}}
\put(82,183){\circle*{8}}
\put(88,180){\scriptsize{-2}}
\put(90,50){\line(1,0){10}}
\put(108,50){\circle{16}}
\put(106,47){\footnotesize{4}}
\put(116,50){\line(1,0){8}}
\put(128,50){\circle*{8}}
\put(125,56){\scriptsize{-2}}
\put(132,50){\line(1,0){8}}
\put(144,50){\circle*{8}}
\put(141,56){\scriptsize{-2}}
\put(148,50){\line(1,0){8}}
\put(160,50){\circle*{8}}
\put(157,56){\scriptsize{-3}}
\put(164,50){\line(1,0){8}}
\put(180,50){\circle{16}}
\put(178,47){\footnotesize{5}}
\put(188,50){\line(1,0){8}}
\put(200,50){\circle*{8}}
\put(197,56){\scriptsize{-4}}
\put(200,46){\line(0,-1){8}}
\put(200,34){\circle*{8}}
\put(206,31){\scriptsize{-2}}
\put(200,30){\line(0,-1){8}}
\put(200,14){\circle{16}}
\put(198,11){\footnotesize{6}}
\put(200,6){\line(0,-1){8}}
\put(200,-6){\circle*{8}}
\put(206,-9){\scriptsize{-7}}
\put(196,-6){\line(-1,0){8}}
\put(180,-6){\circle{16}}
\put(178,-9){\footnotesize{7}}
\put(172,-6){\line(-1,0){8}}
\put(160,-6){\circle*{8}}
\put(157,0){\scriptsize{-2}}
\put(156,-6){\line(-1,0){8}}
\put(144,-6){\circle*{8}}
\put(141,0){\scriptsize{-4}}
\put(140,-6){\line(-1,0){8}}
\put(124,-6){\circle{16}}
\put(122,-9){\footnotesize{8}}
\put(116,-6){\line(-1,0){8}}
\put(104,-6){\circle*{8}}
\put(101,0){\scriptsize{-3}}
\put(100,-6){\line(-1,0){8}}
\put(88,-6){\circle*{8}}
\put(85,0){\scriptsize{-2}}
\put(84,-6){\line(-1,0){8}}
\put(72,-6){\circle*{8}}
\put(69,0){\scriptsize{-2}}
\put(68,-6){\line(-1,0){8}}
\put(52,-6){\circle{16}}
\put(50,-9){\footnotesize{9}}
\put(52,2){\line(1,2){22}}
\end{picture}
\end{center}
\vspace{0.2cm}
\begin{proof} In the figure it is shown the extraction of all curves with
discrepancy
0 for the pair
$(\overline{S},\frac67C)$.
The required curves are enumerated and they are
$(-1)$ curves on a minimal resolution. The extraction process was described
in proposition 2.1 \cite{Kud4}.
\par
Let $f\colon \widetilde S\to S$ be an extraction of some set
$\TT$.
The proper transform of $C$ is denoted by $\widetilde{C}$.
The obvious requirement for the set
$\TT$ is
$\widetilde{C}^2<0$. A necessary and sufficient condition for the set
$\TT$ is that the configuration of curves including $\widetilde{C}$ be
contracted to a log terminal singularity.
Therefore all such configurations must have the type
$\AAA_n$, $\DDD_n$, $\EEE_6$, $\EEE_7$ or $\EEE_8$.
Using these facts the reader will easily prove this theorem.
\end{proof}

\begin{corollary} In the case $6-2({\rm ell})$ we have $I=7$, $2\le \rho(S)\le 9$ and
$\rank\Delta(S)+\rho(S)=10$.
\begin{proof} It is proved as theorem 2.8 \cite{Kud4}.
\end{proof}
\end{corollary}

\begin{remark}
The remained statements are proved similarly. The proofs are left to the reader.
\end{remark}

\begin{theorem} In the case $8-1({\rm ell})$ the set $\TT$
must be one of the following ones.
\begin{enumerate}
\item Let $T_1=0$. Then either $3\in \TT_2$, or $\{4,7\}\subset \TT_2$.
\item Let $1\in \TT_1$. Then $\TT_2\ne \emptyset$.
\item Let $\TT_1=\{2\}$. Then $\{3,4\}\cap \TT_2\ne \emptyset$.
\end{enumerate}
\end{theorem}

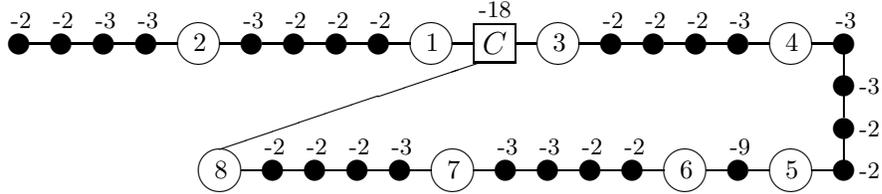
\begin{center}
\begin{picture}(310,100)(0,0)
\put(4,90){\circle*{8}}
\put(1,96){\scriptsize{-3}}
\put(4,50){\circle*{8}}
\put(1,56){\scriptsize{-2}}
\put(8,50){\line(1,0){8}}
\put(20,50){\circle*{8}}
\put(17,56){\scriptsize{-2}}
\put(24,50){\line(1,0){8}}
\put(36,50){\circle*{8}}
\put(33,56){\scriptsize{-3}}
\put(40,50){\line(1,0){8}}
\put(52,50){\circle*{8}}
\put(49,56){\scriptsize{-3}}
\put(56,50){\line(1,0){8}}
\put(72,50){\circle{16}}
\put(70,47){\footnotesize{2}}
\put(80,50){\line(1,0){8}}
\put(92,50){\circle*{8}}
\put(89,56){\scriptsize{-3}}
\put(96,50){\line(1,0){8}}
\put(108,50){\circle*{8}}
\put(105,56){\scriptsize{-2}}
\put(112,50){\line(1,0){8}}
\put(124,50){\circle*{8}}
\put(121,56){\scriptsize{-2}}
\put(128,50){\line(1,0){8}}
\put(140,50){\circle*{8}}
\put(137,56){\scriptsize{-2}}
\put(144,50){\line(1,0){8}}
\put(160,50){\circle{16}}
\put(158,47){\footnotesize{1}}
\put(168,50){\line(1,0){8}}
\put(176,46){\fbox{$C$}}
\put(178,60){\scriptsize{-18}}
\put(192,50){\line(1,0){8}}
\put(208,50){\circle{16}}
\put(206,47){\footnotesize{3}}
\put(216,50){\line(1,0){8}}
\put(228,50){\circle*{8}}
\put(225,56){\scriptsize{-2}}
\put(232,50){\line(1,0){8}}
\put(244,50){\circle*{8}}
\put(241,56){\scriptsize{-2}}
\put(248,50){\line(1,0){8}}
\put(260,50){\circle*{8}}
\put(257,56){\scriptsize{-2}}
\put(264,50){\line(1,0){8}}
\put(276,50){\circle*{8}}
\put(273,56){\scriptsize{-3}}
\put(280,50){\line(1,0){8}}
\put(296,50){\circle{16}}
\put(294,47){\footnotesize{4}}
\put(304,50){\line(1,0){8}}
\put(316,50){\circle*{8}}
\put(313,56){\scriptsize{-3}}
\put(316,46){\line(0,-1){8}}
\put(316,34){\circle*{8}}
\put(322,31){\scriptsize{-3}}
\put(316,30){\line(0,-1){8}}
\put(316,18){\circle*{8}}
\put(322,15){\scriptsize{-2}}
\put(316,14){\line(0,-1){8}}
\put(316,2){\circle*{8}}
\put(322,-1){\scriptsize{-2}}
\put(312,2){\line(-1,0){8}}
\put(296,2){\circle{16}}
\put(294,-1){\footnotesize{5}}
\put(288,2){\line(-1,0){8}}
\put(276,2){\circle*{8}}
\put(273,8){\scriptsize{-9}}
\put(272,2){\line(-1,0){8}}
\put(256,2){\circle{16}}
\put(254,-1){\footnotesize{6}}
\put(248,2){\line(-1,0){8}}
\put(236,2){\circle*{8}}
\put(233,8){\scriptsize{-2}}
\put(232,2){\line(-1,0){8}}
\put(220,2){\circle*{8}}
\put(217,8){\scriptsize{-2}}
\put(216,2){\line(-1,0){8}}
\put(204,2){\circle*{8}}
\put(201,8){\scriptsize{-3}}
\put(200,2){\line(-1,0){8}}
\put(188,2){\circle*{8}}
\put(185,8){\scriptsize{-3}}
\put(184,2){\line(-1,0){8}}
\put(168,2){\circle{16}}
\put(166,-1){\footnotesize{7}}
\put(160,2){\line(-1,0){8}}
\put(148,2){\circle*{8}}
\put(145,8){\scriptsize{-3}}
\put(144,2){\line(-1,0){8}}
\put(132,2){\circle*{8}}
\put(129,8){\scriptsize{-2}}
\put(128,2){\line(-1,0){8}}
\put(116,2){\circle*{8}}
\put(113,8){\scriptsize{-2}}
\put(112,2){\line(-1,0){8}}
\put(100,2){\circle*{8}}
\put(97,8){\scriptsize{-2}}
\put(96,2){\line(-1,0){8}}
\put(80,2){\circle{16}}
\put(78,-1){\footnotesize{8}}
\put(80,10){\line(3,1){98}}
\end{picture}
\end{center}

\vspace{0.2cm}
\begin{corollary} In the case $8-1({\rm ell})$ we have $I=9$, $1\le \rho(S)\le 8$ and
$\rank\Delta(S)+\rho(S)=12$.
\end{corollary}

\begin{theorem} In the case $9-1({\rm ell})$ the set $\TT$
must be one of the following ones.
\begin{enumerate}
\item Let $T_1=0$. Then $2\in \TT_2$.
\item Let $\TT_1=\{1\}$. Then $\TT_2\ne \emptyset$.
\end{enumerate}
\end{theorem}

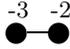
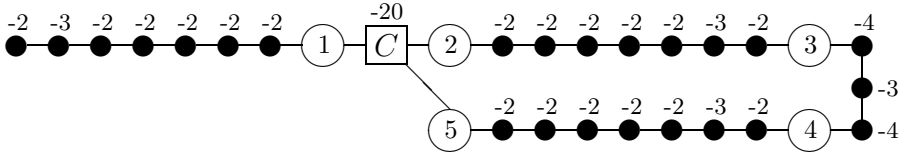
\begin{center}
\begin{picture}(330,100)(0,0)
\put(4,90){\circle*{8}}
\put(1,96){\scriptsize{-3}}
\put(8,90){\line(1,0){8}}
\put(20,90){\circle*{8}}
\put(17,96){\scriptsize{-2}}
\put(4,50){\circle*{8}}
\put(1,56){\scriptsize{-2}}
\put(8,50){\line(1,0){8}}
\put(20,50){\circle*{8}}
\put(17,56){\scriptsize{-3}}
\put(24,50){\line(1,0){8}}
\put(36,50){\circle*{8}}
\put(33,56){\scriptsize{-2}}
\put(40,50){\line(1,0){8}}
\put(52,50){\circle*{8}}
\put(49,56){\scriptsize{-2}}
\put(56,50){\line(1,0){8}}
\put(68,50){\circle*{8}}
\put(65,56){\scriptsize{-2}}
\put(72,50){\line(1,0){8}}
\put(84,50){\circle*{8}}
\put(81,56){\scriptsize{-2}}
\put(88,50){\line(1,0){8}}
\put(100,50){\circle*{8}}
\put(97,56){\scriptsize{-2}}
\put(104,50){\line(1,0){8}}
\put(120,50){\circle{16}}
\put(118,47){\footnotesize{1}}
\put(128,50){\line(1,0){8}}
\put(136,46){\fbox{$C$}}
\put(138,60){\scriptsize{-20}}
\put(152,50){\line(1,0){8}}
\put(168,50){\circle{16}}
\put(166,47){\footnotesize{2}}
\put(176,50){\line(1,0){8}}
\put(188,50){\circle*{8}}
\put(185,56){\scriptsize{-2}}
\put(192,50){\line(1,0){8}}
\put(204,50){\circle*{8}}
\put(201,56){\scriptsize{-2}}
\put(208,50){\line(1,0){8}}
\put(220,50){\circle*{8}}
\put(217,56){\scriptsize{-2}}
\put(224,50){\line(1,0){8}}
\put(236,50){\circle*{8}}
\put(233,56){\scriptsize{-2}}
\put(240,50){\line(1,0){8}}
\put(252,50){\circle*{8}}
\put(249,56){\scriptsize{-2}}
\put(256,50){\line(1,0){8}}
\put(268,50){\circle*{8}}
\put(265,56){\scriptsize{-3}}
\put(272,50){\line(1,0){8}}
\put(284,50){\circle*{8}}
\put(281,56){\scriptsize{-2}}
\put(288,50){\line(1,0){8}}
\put(304,50){\circle{16}}
\put(302,47){\footnotesize{3}}
\put(312,50){\line(1,0){8}}
\put(324,50){\circle*{8}}
\put(321,56){\scriptsize{-4}}
\put(324,46){\line(0,-1){8}}
\put(324,34){\circle*{8}}
\put(330,31){\scriptsize{-3}}
\put(324,30){\line(0,-1){8}}
\put(324,18){\circle*{8}}
\put(330,15){\scriptsize{-4}}
\put(168,26){\line(-1,1){17}}
\put(168,18){\circle{16}}
\put(166,15){\footnotesize{5}}
\put(176,18){\line(1,0){8}}
\put(188,18){\circle*{8}}
\put(185,24){\scriptsize{-2}}
\put(192,18){\line(1,0){8}}
\put(204,18){\circle*{8}}
\put(201,24){\scriptsize{-2}}
\put(208,18){\line(1,0){8}}
\put(220,18){\circle*{8}}
\put(217,24){\scriptsize{-2}}
\put(224,18){\line(1,0){8}}
\put(236,18){\circle*{8}}
\put(233,24){\scriptsize{-2}}
\put(240,18){\line(1,0){8}}
\put(252,18){\circle*{8}}
\put(249,24){\scriptsize{-2}}
\put(256,18){\line(1,0){8}}
\put(268,18){\circle*{8}}
\put(265,24){\scriptsize{-3}}
\put(272,18){\line(1,0){8}}
\put(284,18){\circle*{8}}
\put(281,24){\scriptsize{-2}}
\put(288,18){\line(1,0){8}}
\put(304,18){\circle{16}}
\put(302,15){\footnotesize{4}}
\put(312,18){\line(1,0){8}}
\end{picture}
\end{center}

\vspace{0.2cm}
\begin{corollary} In the case $9-1({\rm ell})$ we have $I=10$, $1\le \rho(S)\le 5$ and
$\rank\Delta(S)+\rho(S)=12$.
\end{corollary}

\begin{theorem} In the case $22-1({\rm ell})$ the set $\TT$
must satisfy the following condition: either
$3\in \TT_2$, or $1\in \TT_1$ and $\TT_2\ne \emptyset$.
\end{theorem}

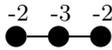
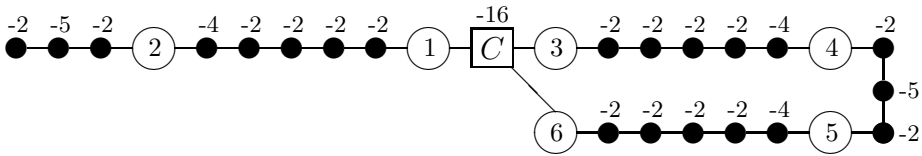
\begin{center}
\begin{picture}(310,100)(10,0)
\put(4,90){\circle*{8}}
\put(1,96){\scriptsize{-2}}
\put(8,90){\line(1,0){8}}
\put(20,90){\circle*{8}}
\put(17,96){\scriptsize{-3}}
\put(24,90){\line(1,0){8}}
\put(36,90){\circle*{8}}
\put(33,96){\scriptsize{-2}}
\put(4,50){\circle*{8}}
\put(1,56){\scriptsize{-2}}
\put(8,50){\line(1,0){8}}
\put(20,50){\circle*{8}}
\put(17,56){\scriptsize{-5}}
\put(24,50){\line(1,0){8}}
\put(36,50){\circle*{8}}
\put(33,56){\scriptsize{-2}}
\put(40,50){\line(1,0){8}}
\put(56,50){\circle{16}}
\put(54,47){\footnotesize{2}}
\put(64,50){\line(1,0){8}}
\put(76,50){\circle*{8}}
\put(73,56){\scriptsize{-4}}
\put(80,50){\line(1,0){8}}
\put(92,50){\circle*{8}}
\put(89,56){\scriptsize{-2}}
\put(96,50){\line(1,0){8}}
\put(108,50){\circle*{8}}
\put(105,56){\scriptsize{-2}}
\put(112,50){\line(1,0){8}}
\put(124,50){\circle*{8}}
\put(121,56){\scriptsize{-2}}
\put(128,50){\line(1,0){8}}
\put(140,50){\circle*{8}}
\put(137,56){\scriptsize{-2}}
\put(144,50){\line(1,0){8}}
\put(160,50){\circle{16}}
\put(158,47){\footnotesize{1}}
\put(168,50){\line(1,0){8}}
\put(176,46){\fbox{$C$}}
\put(178,60){\scriptsize{-16}}
\put(192,50){\line(1,0){8}}
\put(208,50){\circle{16}}
\put(206,47){\footnotesize{3}}
\put(216,50){\line(1,0){8}}
\put(228,50){\circle*{8}}
\put(225,56){\scriptsize{-2}}
\put(232,50){\line(1,0){8}}
\put(244,50){\circle*{8}}
\put(241,56){\scriptsize{-2}}
\put(248,50){\line(1,0){8}}
\put(260,50){\circle*{8}}
\put(257,56){\scriptsize{-2}}
\put(264,50){\line(1,0){8}}
\put(276,50){\circle*{8}}
\put(273,56){\scriptsize{-2}}
\put(280,50){\line(1,0){8}}
\put(292,50){\circle*{8}}
\put(289,56){\scriptsize{-4}}
\put(296,50){\line(1,0){8}}
\put(312,50){\circle{16}}
\put(310,47){\footnotesize{4}}
\put(320,50){\line(1,0){8}}
\put(332,50){\circle*{8}}
\put(329,56){\scriptsize{-2}}
\put(332,46){\line(0,-1){8}}
\put(332,34){\circle*{8}}
\put(338,31){\scriptsize{-5}}
\put(332,30){\line(0,-1){8}}
\put(332,18){\circle*{8}}
\put(338,15){\scriptsize{-2}}
\put(208,18){\circle{16}}
\put(206,15){\footnotesize{6}}
\put(216,18){\line(1,0){8}}
\put(228,18){\circle*{8}}
\put(225,24){\scriptsize{-2}}
\put(232,18){\line(1,0){8}}
\put(244,18){\circle*{8}}
\put(241,24){\scriptsize{-2}}
\put(248,18){\line(1,0){8}}
\put(260,18){\circle*{8}}
\put(257,24){\scriptsize{-2}}
\put(264,18){\line(1,0){8}}
\put(276,18){\circle*{8}}
\put(273,24){\scriptsize{-2}}
\put(280,18){\line(1,0){8}}
\put(292,18){\circle*{8}}
\put(289,24){\scriptsize{-4}}
\put(296,18){\line(1,0){8}}
\put(312,18){\circle{16}}
\put(310,15){\footnotesize{5}}
\put(320,18){\line(1,0){8}}
\put(332,18){\circle*{8}}
\put(208,26){\line(-1,1){17}}
\end{picture}

\end{center}

\vspace{0.2cm}
\begin{corollary} In the case $22-1({\rm ell})$ we have $I=8$, $1\le \rho(S)\le 6$ and
$\rank\Delta(S)+\rho(S)=12$.
\end{corollary}

\begin{theorem} In the case $51-2({\rm ell})$ the set $\TT$
must be one of the following ones.
\begin{enumerate}
\item Let $T_1=0$. Then either $\{2,3\}\cap\TT_2\ne \emptyset$, or $4\in\TT_2$ and
$\overline{T_2}\ge 7$.
\item Let $\TT_1=\{1\}$. Then $\TT_2\ne \emptyset$.
\end{enumerate}
\end{theorem}

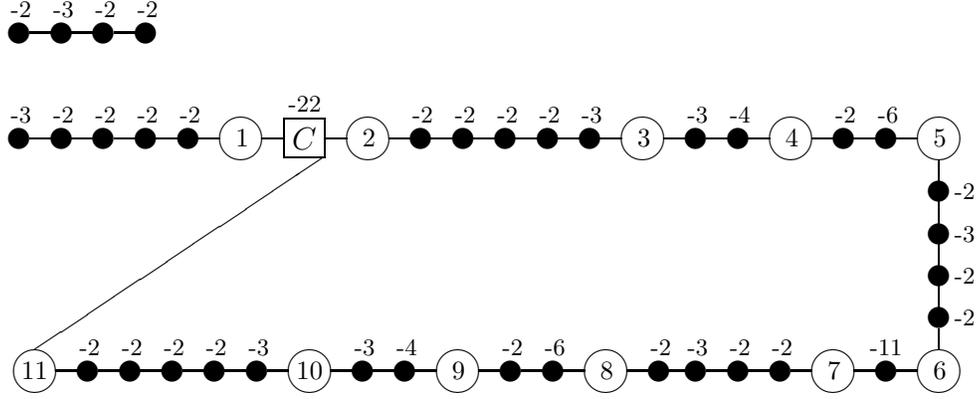
\begin{center}
\begin{picture}(360,150)(0,0)
\put(4,140){\circle*{8}}
\put(1,146){\scriptsize{-2}}
\put(8,140){\line(1,0){8}}
\put(20,140){\circle*{8}}
\put(17,146){\scriptsize{-3}}
\put(24,140){\line(1,0){8}}
\put(36,140){\circle*{8}}
\put(33,146){\scriptsize{-2}}
\put(40,140){\line(1,0){8}}
\put(52,140){\circle*{8}}
\put(49,146){\scriptsize{-2}}
\put(4,100){\circle*{8}}
\put(1,106){\scriptsize{-3}}
\put(8,100){\line(1,0){8}}
\put(20,100){\circle*{8}}
\put(17,106){\scriptsize{-2}}
\put(24,100){\line(1,0){8}}
\put(36,100){\circle*{8}}
\put(33,106){\scriptsize{-2}}
\put(40,100){\line(1,0){8}}
\put(52,100){\circle*{8}}
\put(49,106){\scriptsize{-2}}
\put(56,100){\line(1,0){8}}
\put(68,100){\circle*{8}}
\put(65,106){\scriptsize{-2}}
\put(72,100){\line(1,0){8}}
\put(88,100){\circle{16}}
\put(86,97){\footnotesize{1}}
\put(96,100){\line(1,0){8}}
\put(104,96){\fbox{$C$}}
\put(106,110){\scriptsize{-22}}
\put(120,100){\line(1,0){8}}
\put(136,100){\circle{16}}
\put(134,97){\footnotesize{2}}
\put(144,100){\line(1,0){8}}
\put(156,100){\circle*{8}}
\put(153,106){\scriptsize{-2}}
\put(160,100){\line(1,0){8}}
\put(172,100){\circle*{8}}
\put(169,106){\scriptsize{-2}}
\put(176,100){\line(1,0){8}}
\put(188,100){\circle*{8}}
\put(185,106){\scriptsize{-2}}
\put(192,100){\line(1,0){8}}
\put(204,100){\circle*{8}}
\put(201,106){\scriptsize{-2}}
\put(208,100){\line(1,0){8}}
\put(220,100){\circle*{8}}
\put(217,106){\scriptsize{-3}}
\put(224,100){\line(1,0){8}}
\put(240,100){\circle{16}}
\put(238,97){\footnotesize{3}}
\put(248,100){\line(1,0){8}}
\put(260,100){\circle*{8}}
\put(257,106){\scriptsize{-3}}
\put(264,100){\line(1,0){8}}
\put(276,100){\circle*{8}}
\put(273,106){\scriptsize{-4}}
\put(280,100){\line(1,0){8}}
\put(296,100){\circle{16}}
\put(294,97){\footnotesize{4}}
\put(304,100){\line(1,0){8}}
\put(316,100){\circle*{8}}
\put(313,106){\scriptsize{-2}}
\put(320,100){\line(1,0){8}}
\put(332,100){\circle*{8}}
\put(329,106){\scriptsize{-6}}
\put(336,100){\line(1,0){8}}
\put(352,100){\circle{16}}
\put(350,97){\footnotesize{5}}
\put(352,92){\line(0,-1){8}}
\put(352,80){\circle*{8}}
\put(358,77){\scriptsize{-2}}
\put(352,76){\line(0,-1){8}}
\put(352,64){\circle*{8}}
\put(358,61){\scriptsize{-3}}
\put(352,60){\line(0,-1){8}}
\put(352,48){\circle*{8}}
\put(358,45){\scriptsize{-2}}
\put(352,44){\line(0,-1){8}}
\put(352,32){\circle*{8}}
\put(358,29){\scriptsize{-2}}
\put(352,28){\line(0,-1){8}}
\put(352,12){\circle{16}}
\put(350,9){\footnotesize{6}}
\put(344,12){\line(-1,0){8}}
\put(332,12){\circle*{8}}
\put(326,18){\scriptsize{-11}}
\put(328,12){\line(-1,0){8}}
\put(312,12){\circle{16}}
\put(310,9){\footnotesize{7}}
\put(304,12){\line(-1,0){8}}
\put(292,12){\circle*{8}}
\put(289,18){\scriptsize{-2}}
\put(288,12){\line(-1,0){8}}
\put(276,12){\circle*{8}}
\put(273,18){\scriptsize{-2}}
\put(272,12){\line(-1,0){8}}
\put(260,12){\circle*{8}}
\put(257,18){\scriptsize{-3}}
\put(256,12){\line(-1,0){8}}
\put(246,12){\circle*{8}}
\put(243,18){\scriptsize{-2}}
\put(242,12){\line(-1,0){8}}
\put(226,12){\circle{16}}
\put(224,9){\footnotesize{8}}
\put(218,12){\line(-1,0){8}}
\put(206,12){\circle*{8}}
\put(203,18){\scriptsize{-6}}
\put(202,12){\line(-1,0){8}}
\put(190,12){\circle*{8}}
\put(187,18){\scriptsize{-2}}
\put(186,12){\line(-1,0){8}}
\put(170,12){\circle{16}}
\put(168,9){\footnotesize{9}}
\put(162,12){\line(-1,0){8}}
\put(150,12){\circle*{8}}
\put(147,18){\scriptsize{-4}}
\put(146,12){\line(-1,0){8}}
\put(134,12){\circle*{8}}
\put(131,18){\scriptsize{-3}}
\put(130,12){\line(-1,0){8}}
\put(114,12){\circle{16}}
\put(109,9){\footnotesize{10}}
\put(106,12){\line(-1,0){8}}
\put(94,12){\circle*{8}}
\put(91,18){\scriptsize{-3}}
\put(90,12){\line(-1,0){8}}
\put(78,12){\circle*{8}}
\put(75,18){\scriptsize{-2}}
\put(74,12){\line(-1,0){8}}
\put(62,12){\circle*{8}}
\put(59,18){\scriptsize{-2}}
\put(58,12){\line(-1,0){8}}
\put(46,12){\circle*{8}}
\put(43,18){\scriptsize{-2}}
\put(42,12){\line(-1,0){8}}
\put(30,12){\circle*{8}}
\put(27,18){\scriptsize{-2}}
\put(26,12){\line(-1,0){8}}
\put(10,12){\circle{16}}
\put(5,9){\footnotesize{11}}
\put(10,20){\line(3,2){109}}
\end{picture}
\end{center}

\vspace{0.2cm}
\begin{corollary} In the case $51-2({\rm ell})$ we have $I=11$, $1\le \rho(S)\le 11$
and $\rank\Delta(S)+\rho(S)=12$.
\end{corollary}

\begin{remark} The $K3$ surfaces with the automorphisms of order 11 and hence
the log Enriques surfaces with index
$I=11$ were described in \cite{OZ2}.
\end{remark}

\begin{theorem} In the case $51-6({\rm ell})$ the set $\TT$
must satisfy the following condition: either
$2\in \TT_2$, or $\TT_1=\{1\}$ and $\TT_2\ne \emptyset$.
\end{theorem}

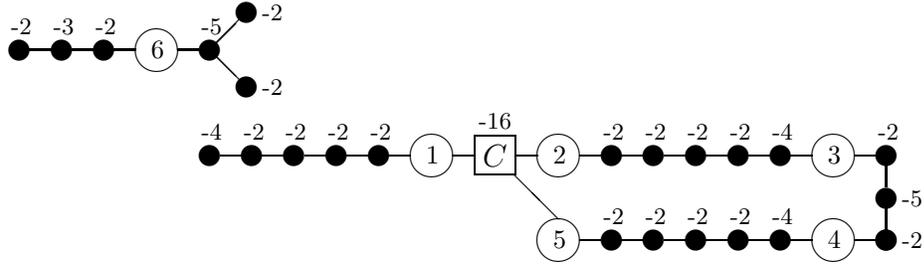
\begin{center}
\begin{picture}(300,120)(0,0)
\put(4,90){\circle*{8}}
\put(1,96){\scriptsize{-2}}
\put(8,90){\line(1,0){8}}
\put(20,90){\circle*{8}}
\put(17,96){\scriptsize{-3}}
\put(24,90){\line(1,0){8}}
\put(36,90){\circle*{8}}
\put(33,96){\scriptsize{-2}}
\put(40,90){\line(1,0){8}}
\put(56,90){\circle{16}}
\put(54,87){\footnotesize{6}}
\put(64,90){\line(1,0){8}}
\put(76,90){\circle*{8}}
\put(73,96){\scriptsize{-5}}
\put(76,90){\line(1,1){12}}
\put(90,104){\circle*{8}}
\put(96,101){\scriptsize{-2}}
\put(76,90){\line(1,-1){12}}
\put(90,76){\circle*{8}}
\put(96,73){\scriptsize{-2}}
\put(76,50){\circle*{8}}
\put(73,56){\scriptsize{-4}}
\put(80,50){\line(1,0){8}}
\put(92,50){\circle*{8}}
\put(89,56){\scriptsize{-2}}
\put(96,50){\line(1,0){8}}
\put(108,50){\circle*{8}}
\put(105,56){\scriptsize{-2}}
\put(112,50){\line(1,0){8}}
\put(124,50){\circle*{8}}
\put(121,56){\scriptsize{-2}}
\put(128,50){\line(1,0){8}}
\put(140,50){\circle*{8}}
\put(137,56){\scriptsize{-2}}
\put(144,50){\line(1,0){8}}
\put(160,50){\circle{16}}
\put(158,47){\footnotesize{1}}
\put(168,50){\line(1,0){8}}
\put(176,46){\fbox{$C$}}
\put(178,60){\scriptsize{-16}}
\put(192,50){\line(1,0){8}}
\put(208,50){\circle{16}}
\put(206,47){\footnotesize{2}}
\put(216,50){\line(1,0){8}}
\put(228,50){\circle*{8}}
\put(225,56){\scriptsize{-2}}
\put(232,50){\line(1,0){8}}
\put(244,50){\circle*{8}}
\put(241,56){\scriptsize{-2}}
\put(248,50){\line(1,0){8}}
\put(260,50){\circle*{8}}
\put(257,56){\scriptsize{-2}}
\put(264,50){\line(1,0){8}}
\put(276,50){\circle*{8}}
\put(273,56){\scriptsize{-2}}
\put(280,50){\line(1,0){8}}
\put(292,50){\circle*{8}}
\put(289,56){\scriptsize{-4}}
\put(296,50){\line(1,0){8}}
\put(312,50){\circle{16}}
\put(310,47){\footnotesize{3}}
\put(320,50){\line(1,0){8}}
\put(332,50){\circle*{8}}
\put(329,56){\scriptsize{-2}}
\put(332,46){\line(0,-1){8}}
\put(332,34){\circle*{8}}
\put(338,31){\scriptsize{-5}}
\put(332,30){\line(0,-1){8}}
\put(332,18){\circle*{8}}
\put(338,15){\scriptsize{-2}}
\put(208,18){\circle{16}}
\put(206,15){\footnotesize{5}}
\put(216,18){\line(1,0){8}}
\put(228,18){\circle*{8}}
\put(225,24){\scriptsize{-2}}
\put(232,18){\line(1,0){8}}
\put(244,18){\circle*{8}}
\put(241,24){\scriptsize{-2}}
\put(248,18){\line(1,0){8}}
\put(260,18){\circle*{8}}
\put(257,24){\scriptsize{-2}}
\put(264,18){\line(1,0){8}}
\put(276,18){\circle*{8}}
\put(273,24){\scriptsize{-2}}
\put(280,18){\line(1,0){8}}
\put(292,18){\circle*{8}}
\put(289,24){\scriptsize{-4}}
\put(296,18){\line(1,0){8}}
\put(312,18){\circle{16}}
\put(310,15){\footnotesize{4}}
\put(320,18){\line(1,0){8}}
\put(332,18){\circle*{8}}
\put(208,26){\line(-1,1){17}}
\end{picture}
\end{center}

\vspace{0.2cm}
\begin{corollary} In the case $51-6({\rm ell})$ we have $I=8$, $1\le \rho(S)\le 6$ and
$\rank\Delta(S)+\rho(S)=12$.
\end{corollary}

\begin{theorem} In the case $52-2({\rm ell})$ the set $\TT$
must satisfy the following condition: either
$\{3,4\}\cap\TT_2\ne \emptyset$, or $T_1\ge 1$ and $\TT_2\ne \emptyset$.
\end{theorem}

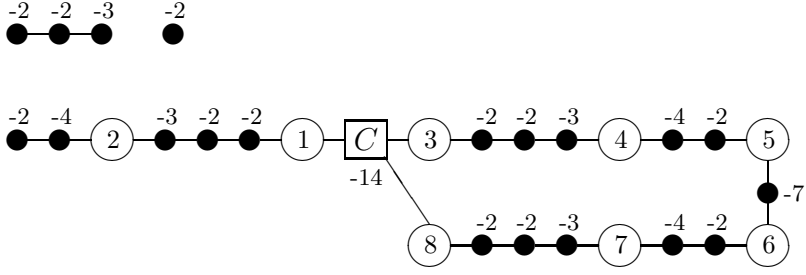
\begin{center}
\begin{picture}(300,100)(0,0)
\put(4,90){\circle*{8}}
\put(1,96){\scriptsize{-2}}
\put(8,90){\line(1,0){8}}
\put(20,90){\circle*{8}}
\put(17,96){\scriptsize{-2}}
\put(24,90){\line(1,0){8}}
\put(36,90){\circle*{8}}
\put(33,96){\scriptsize{-3}}
\put(63,90){\circle*{8}}
\put(60,96){\scriptsize{-2}}
\put(4,50){\circle*{8}}
\put(1,56){\scriptsize{-2}}
\put(8,50){\line(1,0){8}}
\put(20,50){\circle*{8}}
\put(17,56){\scriptsize{-4}}
\put(24,50){\line(1,0){8}}
\put(40,50){\circle{16}}
\put(38,47){\footnotesize{2}}
\put(48,50){\line(1,0){8}}
\put(60,50){\circle*{8}}
\put(57,56){\scriptsize{-3}}
\put(64,50){\line(1,0){8}}
\put(76,50){\circle*{8}}
\put(73,56){\scriptsize{-2}}
\put(80,50){\line(1,0){8}}
\put(92,50){\circle*{8}}
\put(89,56){\scriptsize{-2}}
\put(96,50){\line(1,0){8}}
\put(112,50){\circle{16}}
\put(110,47){\footnotesize{1}}
\put(120,50){\line(1,0){8}}
\put(128,46){\fbox{$C$}}
\put(130,33){\scriptsize{-14}}
\put(144,50){\line(1,0){8}}
\put(160,50){\circle{16}}
\put(158,47){\footnotesize{3}}
\put(168,50){\line(1,0){8}}
\put(180,50){\circle*{8}}
\put(178,56){\scriptsize{-2}}
\put(184,50){\line(1,0){8}}
\put(196,50){\circle*{8}}
\put(193,56){\scriptsize{-2}}
\put(200,50){\line(1,0){8}}
\put(212,50){\circle*{8}}
\put(209,56){\scriptsize{-3}}
\put(216,50){\line(1,0){8}}
\put(232,50){\circle{16}}
\put(230,47){\footnotesize{4}}
\put(240,50){\line(1,0){8}}
\put(252,50){\circle*{8}}
\put(249,56){\scriptsize{-4}}
\put(256,50){\line(1,0){8}}
\put(268,50){\circle*{8}}
\put(265,56){\scriptsize{-2}}
\put(272,50){\line(1,0){8}}
\put(288,50){\circle{16}}
\put(286,47){\footnotesize{5}}
\put(288,42){\line(0,-1){8}}
\put(288,30){\circle*{8}}
\put(294,27){\scriptsize{-7}}
\put(288,26){\line(0,-1){8}}
\put(288,10){\circle{16}}
\put(286,7){\footnotesize{6}}
\put(216,10){\line(1,0){8}}
\put(232,10){\circle{16}}
\put(230,7){\footnotesize{7}}
\put(240,10){\line(1,0){8}}
\put(252,10){\circle*{8}}
\put(249,16){\scriptsize{-4}}
\put(256,10){\line(1,0){8}}
\put(268,10){\circle*{8}}
\put(265,16){\scriptsize{-2}}
\put(272,10){\line(1,0){8}}
\put(160,18){\line(-2,3){17}}
\put(160,10){\circle{16}}
\put(158,7){\footnotesize{8}}
\put(168,10){\line(1,0){8}}
\put(180,10){\circle*{8}}
\put(178,16){\scriptsize{-2}}
\put(184,10){\line(1,0){8}}
\put(196,10){\circle*{8}}
\put(193,16){\scriptsize{-2}}
\put(200,10){\line(1,0){8}}
\put(212,10){\circle*{8}}
\put(209,16){\scriptsize{-3}}
\end{picture}
\end{center}

\vspace{0.2cm}
\begin{corollary} In the case $52-2({\rm ell})$ we have $I=7$, $1\le \rho(S)\le 8$ and
$\rank\Delta(S)+\rho(S)=10$.
\end{corollary}

\begin{theorem} In the case $53-2({\rm ell})$ the set $\TT$
must be one of the following ones.
\begin{enumerate}
\item Let $T_1=0$. Then $\{2,3\}\cap\TT_2\ne \emptyset$.
\item Let $\TT_1=\{1\}$. Then $\TT_2\ne \emptyset$.
\end{enumerate}
\end{theorem}

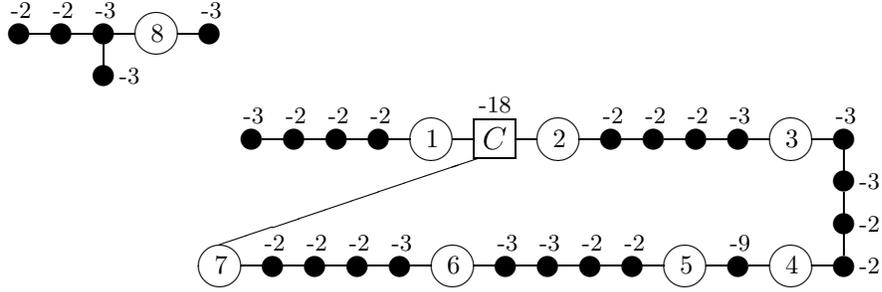
\begin{center}
\begin{picture}(310,100)(0,0)
\put(4,90){\circle*{8}}
\put(1,96){\scriptsize{-2}}
\put(8,90){\line(1,0){8}}
\put(20,90){\circle*{8}}
\put(17,96){\scriptsize{-2}}
\put(24,90){\line(1,0){8}}
\put(36,90){\circle*{8}}
\put(33,96){\scriptsize{-3}}
\put(40,90){\line(1,0){8}}
\put(56,90){\circle{16}}
\put(54,87){\footnotesize{8}}
\put(64,90){\line(1,0){8}}
\put(76,90){\circle*{8}}
\put(73,96){\scriptsize{-3}}
\put(36,86){\line(0,-1){8}}
\put(36,74){\circle*{8}}
\put(42,71){\scriptsize{-3}}
\put(92,50){\circle*{8}}
\put(89,56){\scriptsize{-3}}
\put(96,50){\line(1,0){8}}
\put(108,50){\circle*{8}}
\put(105,56){\scriptsize{-2}}
\put(112,50){\line(1,0){8}}
\put(124,50){\circle*{8}}
\put(121,56){\scriptsize{-2}}
\put(128,50){\line(1,0){8}}
\put(140,50){\circle*{8}}
\put(137,56){\scriptsize{-2}}
\put(144,50){\line(1,0){8}}
\put(160,50){\circle{16}}
\put(158,47){\footnotesize{1}}
\put(168,50){\line(1,0){8}}
\put(176,46){\fbox{$C$}}
\put(178,60){\scriptsize{-18}}
\put(192,50){\line(1,0){8}}
\put(208,50){\circle{16}}
\put(206,47){\footnotesize{2}}
\put(216,50){\line(1,0){8}}
\put(228,50){\circle*{8}}
\put(225,56){\scriptsize{-2}}
\put(232,50){\line(1,0){8}}
\put(244,50){\circle*{8}}
\put(241,56){\scriptsize{-2}}
\put(248,50){\line(1,0){8}}
\put(260,50){\circle*{8}}
\put(257,56){\scriptsize{-2}}
\put(264,50){\line(1,0){8}}
\put(276,50){\circle*{8}}
\put(273,56){\scriptsize{-3}}
\put(280,50){\line(1,0){8}}
\put(296,50){\circle{16}}
\put(294,47){\footnotesize{3}}
\put(304,50){\line(1,0){8}}
\put(316,50){\circle*{8}}
\put(313,56){\scriptsize{-3}}
\put(316,46){\line(0,-1){8}}
\put(316,34){\circle*{8}}
\put(322,31){\scriptsize{-3}}
\put(316,30){\line(0,-1){8}}
\put(316,18){\circle*{8}}
\put(322,15){\scriptsize{-2}}
\put(316,14){\line(0,-1){8}}
\put(316,2){\circle*{8}}
\put(322,-1){\scriptsize{-2}}
\put(312,2){\line(-1,0){8}}
\put(296,2){\circle{16}}
\put(294,-1){\footnotesize{4}}
\put(288,2){\line(-1,0){8}}
\put(276,2){\circle*{8}}
\put(273,8){\scriptsize{-9}}
\put(272,2){\line(-1,0){8}}
\put(256,2){\circle{16}}
\put(254,-1){\footnotesize{5}}
\put(248,2){\line(-1,0){8}}
\put(236,2){\circle*{8}}
\put(233,8){\scriptsize{-2}}
\put(232,2){\line(-1,0){8}}
\put(220,2){\circle*{8}}
\put(217,8){\scriptsize{-2}}
\put(216,2){\line(-1,0){8}}
\put(204,2){\circle*{8}}
\put(201,8){\scriptsize{-3}}
\put(200,2){\line(-1,0){8}}
\put(188,2){\circle*{8}}
\put(185,8){\scriptsize{-3}}
\put(184,2){\line(-1,0){8}}
\put(168,2){\circle{16}}
\put(166,-1){\footnotesize{6}}
\put(160,2){\line(-1,0){8}}
\put(148,2){\circle*{8}}
\put(145,8){\scriptsize{-3}}
\put(144,2){\line(-1,0){8}}
\put(132,2){\circle*{8}}
\put(129,8){\scriptsize{-2}}
\put(128,2){\line(-1,0){8}}
\put(116,2){\circle*{8}}
\put(113,8){\scriptsize{-2}}
\put(112,2){\line(-1,0){8}}
\put(100,2){\circle*{8}}
\put(97,8){\scriptsize{-2}}
\put(96,2){\line(-1,0){8}}
\put(80,2){\circle{16}}
\put(78,-1){\footnotesize{7}}
\put(80,10){\line(3,1){98}}
\end{picture}
\end{center}

\vspace{0.2cm}
\begin{corollary} In the case $53-2({\rm ell})$ we have $I=9$, $1\le \rho(S)\le 8$ and
$\rank\Delta(S)+\rho(S)=12$.
\end{corollary}

\begin{theorem} In the case $54({\rm ell})$ the set $\TT$
must satisfy the following condition: $\TT_2\ne\emptyset$.
\end{theorem}

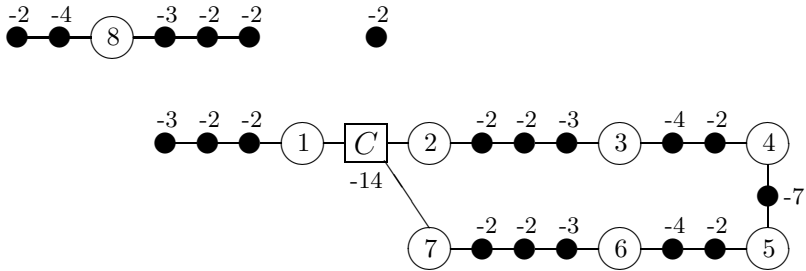
\begin{center}
\begin{picture}(300,100)(0,0)
\put(4,90){\circle*{8}}
\put(1,96){\scriptsize{-2}}
\put(8,90){\line(1,0){8}}
\put(20,90){\circle*{8}}
\put(17,96){\scriptsize{-4}}
\put(24,90){\line(1,0){8}}
\put(40,90){\circle{16}}
\put(38,87){\footnotesize{8}}
\put(48,90){\line(1,0){8}}
\put(60,90){\circle*{8}}
\put(57,96){\scriptsize{-3}}
\put(64,90){\line(1,0){8}}
\put(76,90){\circle*{8}}
\put(73,96){\scriptsize{-2}}
\put(80,90){\line(1,0){8}}
\put(92,90){\circle*{8}}
\put(89,96){\scriptsize{-2}}
\put(140,90){\circle*{8}}
\put(137,96){\scriptsize{-2}}
\put(60,50){\circle*{8}}
\put(57,56){\scriptsize{-3}}
\put(64,50){\line(1,0){8}}
\put(76,50){\circle*{8}}
\put(73,56){\scriptsize{-2}}
\put(80,50){\line(1,0){8}}
\put(92,50){\circle*{8}}
\put(89,56){\scriptsize{-2}}
\put(96,50){\line(1,0){8}}
\put(112,50){\circle{16}}
\put(110,47){\footnotesize{1}}
\put(120,50){\line(1,0){8}}
\put(128,46){\fbox{$C$}}
\put(130,33){\scriptsize{-14}}
\put(144,50){\line(1,0){8}}
\put(160,50){\circle{16}}
\put(158,47){\footnotesize{2}}
\put(168,50){\line(1,0){8}}
\put(180,50){\circle*{8}}
\put(178,56){\scriptsize{-2}}
\put(184,50){\line(1,0){8}}
\put(196,50){\circle*{8}}
\put(193,56){\scriptsize{-2}}
\put(200,50){\line(1,0){8}}
\put(212,50){\circle*{8}}
\put(209,56){\scriptsize{-3}}
\put(216,50){\line(1,0){8}}
\put(232,50){\circle{16}}
\put(230,47){\footnotesize{3}}
\put(240,50){\line(1,0){8}}
\put(252,50){\circle*{8}}
\put(249,56){\scriptsize{-4}}
\put(256,50){\line(1,0){8}}
\put(268,50){\circle*{8}}
\put(265,56){\scriptsize{-2}}
\put(272,50){\line(1,0){8}}
\put(288,50){\circle{16}}
\put(286,47){\footnotesize{4}}
\put(288,42){\line(0,-1){8}}
\put(288,30){\circle*{8}}
\put(294,27){\scriptsize{-7}}
\put(288,26){\line(0,-1){8}}
\put(288,10){\circle{16}}
\put(286,7){\footnotesize{5}}
\put(216,10){\line(1,0){8}}
\put(232,10){\circle{16}}
\put(230,7){\footnotesize{6}}
\put(240,10){\line(1,0){8}}
\put(252,10){\circle*{8}}
\put(249,16){\scriptsize{-4}}
\put(256,10){\line(1,0){8}}
\put(268,10){\circle*{8}}
\put(265,16){\scriptsize{-2}}
\put(272,10){\line(1,0){8}}
\put(160,18){\line(-2,3){17}}
\put(160,10){\circle{16}}
\put(158,7){\footnotesize{7}}
\put(168,10){\line(1,0){8}}
\put(180,10){\circle*{8}}
\put(178,16){\scriptsize{-2}}
\put(184,10){\line(1,0){8}}
\put(196,10){\circle*{8}}
\put(193,16){\scriptsize{-2}}
\put(200,10){\line(1,0){8}}
\put(212,10){\circle*{8}}
\put(209,16){\scriptsize{-3}}
\end{picture}
\end{center}

\vspace{0.2cm}
\begin{corollary} In the case $54({\rm ell})$ we have $I=7$, $1\le \rho(S)\le 8$ and
$\rank\Delta(S)+\rho(S)=10$.
\end{corollary}

\begin{definition}
Let $p_a(C)=0$. These are the cases (8)--(10) of theorem \ref{ldmain}.
In all cases the surface $\overline{S}$ has only three singular points lying on
$C$. In the corresponding figures these points  are enumerated from top to bottom
of the figure.
Then $\TT=\TT_1\cup\TT_2\cup\TT_3$, where
$\TT_i$ consists of curves, which contract to the singular point of $\overline{S}$ with
the number $i$.
\end{definition}

\begin{theorem} In the case $18-1(+1)$ the set
$\TT$ must be one of the following ones.
\begin{enumerate}
\item Let $T_1=0$. Then either $\{3,4\}\cap\TT_2\ne \emptyset$, or
$\{6,7\}\cap\TT_3\ne \emptyset$.
\item Let $1\in \TT_1$. Then $\TT_2$ and $\TT_3$ are arbitrary.
\item Let $\TT_1=\{2\}$. Then either $\TT_2\ne \emptyset$, or
$\underline{T_3}\le 8$.
\end{enumerate}
\end{theorem}

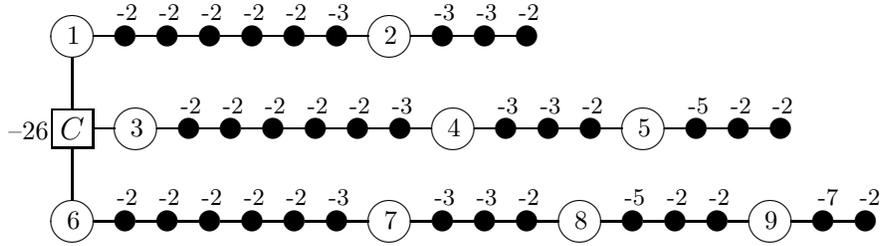
\begin{center}
\begin{picture}(332,90)(0,0)
\put(4,40){\footnotesize{--26}}
\put(20,40){\fbox{$C$}}
\put(28,51){\line(0,1){20}}
\put(28,79){\circle{16}}
\put(26,76){\footnotesize{1}}
\put(36,79){\line(1,0){8}}
\put(48,79){\circle*{8}}
\put(45,85){\scriptsize{-2}}
\put(52,79){\line(1,0){8}}
\put(64,79){\circle*{8}}
\put(61,85){\scriptsize{-2}}
\put(68,79){\line(1,0){8}}
\put(80,79){\circle*{8}}
\put(77,85){\scriptsize{-2}}
\put(84,79){\line(1,0){8}}
\put(96,79){\circle*{8}}
\put(93,85){\scriptsize{-2}}
\put(100,79){\line(1,0){8}}
\put(112,79){\circle*{8}}
\put(109,85){\scriptsize{-2}}
\put(116,79){\line(1,0){8}}
\put(128,79){\circle*{8}}
\put(125,85){\scriptsize{-3}}
\put(132,79){\line(1,0){8}}
\put(148,79){\circle{16}}
\put(146,76){\footnotesize{2}}
\put(156,79){\line(1,0){8}}
\put(168,79){\circle*{8}}
\put(165,85){\scriptsize{-3}}
\put(172,79){\line(1,0){8}}
\put(184,79){\circle*{8}}
\put(181,85){\scriptsize{-3}}
\put(188,79){\line(1,0){8}}
\put(200,79){\circle*{8}}
\put(197,85){\scriptsize{-2}}
\put(36,44){\line(1,0){8}}
\put(52,44){\circle{16}}
\put(50,41){\footnotesize{3}}
\put(60,44){\line(1,0){8}}
\put(72,44){\circle*{8}}
\put(69,50){\scriptsize{-2}}
\put(76,44){\line(1,0){8}}
\put(88,44){\circle*{8}}
\put(85,50){\scriptsize{-2}}
\put(92,44){\line(1,0){8}}
\put(104,44){\circle*{8}}
\put(101,50){\scriptsize{-2}}
\put(108,44){\line(1,0){8}}
\put(120,44){\circle*{8}}
\put(117,50){\scriptsize{-2}}
\put(124,44){\line(1,0){8}}
\put(136,44){\circle*{8}}
\put(133,50){\scriptsize{-2}}
\put(140,44){\line(1,0){8}}
\put(152,44){\circle*{8}}
\put(149,50){\scriptsize{-3}}
\put(156,44){\line(1,0){8}}
\put(172,44){\circle{16}}
\put(170,41){\footnotesize{4}}
\put(180,44){\line(1,0){8}}
\put(192,44){\circle*{8}}
\put(189,50){\scriptsize{-3}}
\put(196,44){\line(1,0){8}}
\put(208,44){\circle*{8}}
\put(205,50){\scriptsize{-3}}
\put(212,44){\line(1,0){8}}
\put(224,44){\circle*{8}}
\put(221,50){\scriptsize{-2}}
\put(228,44){\line(1,0){8}}
\put(244,44){\circle{16}}
\put(242,41){\footnotesize{5}}
\put(252,44){\line(1,0){8}}
\put(264,44){\circle*{8}}
\put(261,50){\scriptsize{-5}}
\put(268,44){\line(1,0){8}}
\put(280,44){\circle*{8}}
\put(277,50){\scriptsize{-2}}
\put(284,44){\line(1,0){8}}
\put(296,44){\circle*{8}}
\put(293,50){\scriptsize{-2}}
\put(28,37){\line(0,-1){20}}
\put(28,9){\circle{16}}
\put(26,6){\footnotesize{6}}
\put(36,9){\line(1,0){8}}
\put(48,9){\circle*{8}}
\put(45,15){\scriptsize{-2}}
\put(52,9){\line(1,0){8}}
\put(64,9){\circle*{8}}
\put(61,15){\scriptsize{-2}}
\put(68,9){\line(1,0){8}}
\put(80,9){\circle*{8}}
\put(77,15){\scriptsize{-2}}
\put(84,9){\line(1,0){8}}
\put(96,9){\circle*{8}}
\put(93,15){\scriptsize{-2}}
\put(100,9){\line(1,0){8}}
\put(112,9){\circle*{8}}
\put(109,15){\scriptsize{-2}}
\put(116,9){\line(1,0){8}}
\put(128,9){\circle*{8}}
\put(125,15){\scriptsize{-3}}
\put(132,9){\line(1,0){8}}
\put(148,9){\circle{16}}
\put(146,6){\footnotesize{7}}
\put(156,9){\line(1,0){8}}
\put(168,9){\circle*{8}}
\put(165,15){\scriptsize{-3}}
\put(172,9){\line(1,0){8}}
\put(184,9){\circle*{8}}
\put(181,15){\scriptsize{-3}}
\put(188,9){\line(1,0){8}}
\put(200,9){\circle*{8}}
\put(197,15){\scriptsize{-2}}
\put(204,9){\line(1,0){8}}
\put(220,9){\circle{16}}
\put(218,6){\footnotesize{8}}
\put(228,9){\line(1,0){8}}
\put(240,9){\circle*{8}}
\put(237,15){\scriptsize{-5}}
\put(244,9){\line(1,0){8}}
\put(256,9){\circle*{8}}
\put(253,15){\scriptsize{-2}}
\put(260,9){\line(1,0){8}}
\put(272,9){\circle*{8}}
\put(269,15){\scriptsize{-2}}
\put(276,9){\line(1,0){8}}
\put(292,9){\circle{16}}
\put(290,6){\footnotesize{9}}
\put(300,9){\line(1,0){8}}
\put(312,9){\circle*{8}}
\put(310,15){\scriptsize{-7}}
\put(316,9){\line(1,0){8}}
\put(328,9){\circle*{8}}
\put(326,15){\scriptsize{-2}}
\end{picture}
\end{center}

\vspace{0.2cm}
\begin{corollary} In the case $18-1(+1)$ we have $I=13$, $1\le \rho(S)\le 9$ and
$\rank\Delta(S)+\rho(S)=10$.
\end{corollary}

\begin{theorem} In the case $25-1(+1)$ the set
$\TT$ must be one of the following ones..
\begin{enumerate}
\item Let $T_1=0$. Then either $2\in \TT_2$, or $4\in \TT_3$.
\item Let $\TT_1=\{1\}$. Then $\TT_2$ and $\TT_3$ are arbitrary.
\end{enumerate}
\end{theorem}

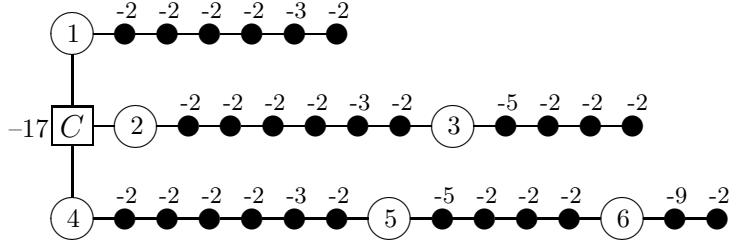
\begin{center}
\begin{picture}(332,90)(0,0)
\put(4,40){\footnotesize{--17}}
\put(20,40){\fbox{$C$}}
\put(28,51){\line(0,1){20}}
\put(28,79){\circle{16}}
\put(26,76){\footnotesize{1}}
\put(36,79){\line(1,0){8}}
\put(48,79){\circle*{8}}
\put(45,85){\scriptsize{-2}}
\put(52,79){\line(1,0){8}}
\put(64,79){\circle*{8}}
\put(61,85){\scriptsize{-2}}
\put(68,79){\line(1,0){8}}
\put(80,79){\circle*{8}}
\put(77,85){\scriptsize{-2}}
\put(84,79){\line(1,0){8}}
\put(96,79){\circle*{8}}
\put(93,85){\scriptsize{-2}}
\put(100,79){\line(1,0){8}}
\put(112,79){\circle*{8}}
\put(109,85){\scriptsize{-3}}
\put(116,79){\line(1,0){8}}
\put(128,79){\circle*{8}}
\put(125,85){\scriptsize{-2}}
\put(36,44){\line(1,0){8}}
\put(52,44){\circle{16}}
\put(50,41){\footnotesize{2}}
\put(60,44){\line(1,0){8}}
\put(72,44){\circle*{8}}
\put(69,50){\scriptsize{-2}}
\put(76,44){\line(1,0){8}}
\put(88,44){\circle*{8}}
\put(85,50){\scriptsize{-2}}
\put(92,44){\line(1,0){8}}
\put(104,44){\circle*{8}}
\put(101,50){\scriptsize{-2}}
\put(108,44){\line(1,0){8}}
\put(120,44){\circle*{8}}
\put(117,50){\scriptsize{-2}}
\put(124,44){\line(1,0){8}}
\put(136,44){\circle*{8}}
\put(133,50){\scriptsize{-3}}
\put(140,44){\line(1,0){8}}
\put(152,44){\circle*{8}}
\put(149,50){\scriptsize{-2}}
\put(156,44){\line(1,0){8}}
\put(172,44){\circle{16}}
\put(170,41){\footnotesize{3}}
\put(180,44){\line(1,0){8}}
\put(192,44){\circle*{8}}
\put(189,50){\scriptsize{-5}}
\put(196,44){\line(1,0){8}}
\put(208,44){\circle*{8}}
\put(205,50){\scriptsize{-2}}
\put(212,44){\line(1,0){8}}
\put(224,44){\circle*{8}}
\put(221,50){\scriptsize{-2}}
\put(228,44){\line(1,0){8}}
\put(240,44){\circle*{8}}
\put(238,50){\scriptsize{-2}}
\put(28,37){\line(0,-1){20}}
\put(28,9){\circle{16}}
\put(26,6){\footnotesize{4}}
\put(36,9){\line(1,0){8}}
\put(48,9){\circle*{8}}
\put(45,15){\scriptsize{-2}}
\put(52,9){\line(1,0){8}}
\put(64,9){\circle*{8}}
\put(61,15){\scriptsize{-2}}
\put(68,9){\line(1,0){8}}
\put(80,9){\circle*{8}}
\put(77,15){\scriptsize{-2}}
\put(84,9){\line(1,0){8}}
\put(96,9){\circle*{8}}
\put(93,15){\scriptsize{-2}}
\put(100,9){\line(1,0){8}}
\put(112,9){\circle*{8}}
\put(109,15){\scriptsize{-3}}
\put(116,9){\line(1,0){8}}
\put(128,9){\circle*{8}}
\put(125,15){\scriptsize{-2}}
\put(132,9){\line(1,0){8}}
\put(148,9){\circle{16}}
\put(146,6){\footnotesize{5}}
\put(156,9){\line(1,0){8}}
\put(168,9){\circle*{8}}
\put(165,15){\scriptsize{-5}}
\put(172,9){\line(1,0){8}}
\put(184,9){\circle*{8}}
\put(181,15){\scriptsize{-2}}
\put(188,9){\line(1,0){8}}
\put(200,9){\circle*{8}}
\put(197,15){\scriptsize{-2}}
\put(204,9){\line(1,0){8}}
\put(216,9){\circle*{8}}
\put(213,15){\scriptsize{-2}}
\put(220,9){\line(1,0){8}}
\put(236,9){\circle{16}}
\put(234,6){\footnotesize{6}}
\put(244,9){\line(1,0){8}}
\put(256,9){\circle*{8}}
\put(253,15){\scriptsize{-9}}
\put(260,9){\line(1,0){8}}
\put(272,9){\circle*{8}}
\put(269,15){\scriptsize{-2}}
\end{picture}
\end{center}

\vspace{0.2cm}
\begin{corollary} In the case $25-1(+1)$ we have $I=17$, $1\le \rho(S)\le 6$ and
$\rank\Delta(S)+\rho(S)=6$.
\end{corollary}

\begin{theorem} In the case $55(0)$ the set
$\TT$
must be one of the following ones.
\begin{enumerate}
\item Let $T_1=0$. Then either $\{2,3\}\cap \TT_2\ne \emptyset$, or
$\underline{T_3}\le 8$, or $4\in\TT_2$ and $\underline{T_3}\le 10$.
\item Let $\TT_1=\{1\}$. Then $\TT_2$ and $\TT_3$ are arbitrary.
\end{enumerate}
\end{theorem}

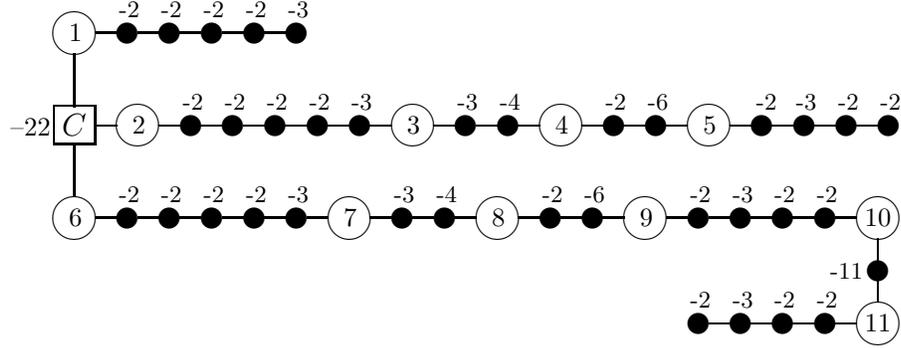
\begin{center}
\begin{picture}(332,120)(0,-30)
\put(4,40){\footnotesize{--22}}
\put(20,40){\fbox{$C$}}
\put(28,51){\line(0,1){20}}
\put(28,79){\circle{16}}
\put(26,76){\footnotesize{1}}
\put(36,79){\line(1,0){8}}
\put(48,79){\circle*{8}}
\put(45,85){\scriptsize{-2}}
\put(52,79){\line(1,0){8}}
\put(64,79){\circle*{8}}
\put(61,85){\scriptsize{-2}}
\put(68,79){\line(1,0){8}}
\put(80,79){\circle*{8}}
\put(77,85){\scriptsize{-2}}
\put(84,79){\line(1,0){8}}
\put(96,79){\circle*{8}}
\put(93,85){\scriptsize{-2}}
\put(100,79){\line(1,0){8}}
\put(112,79){\circle*{8}}
\put(109,85){\scriptsize{-3}}
\put(36,44){\line(1,0){8}}
\put(52,44){\circle{16}}
\put(50,41){\footnotesize{2}}
\put(60,44){\line(1,0){8}}
\put(72,44){\circle*{8}}
\put(69,50){\scriptsize{-2}}
\put(76,44){\line(1,0){8}}
\put(88,44){\circle*{8}}
\put(85,50){\scriptsize{-2}}
\put(92,44){\line(1,0){8}}
\put(104,44){\circle*{8}}
\put(101,50){\scriptsize{-2}}
\put(108,44){\line(1,0){8}}
\put(120,44){\circle*{8}}
\put(117,50){\scriptsize{-2}}
\put(124,44){\line(1,0){8}}
\put(136,44){\circle*{8}}
\put(133,50){\scriptsize{-3}}
\put(140,44){\line(1,0){8}}
\put(156,44){\circle{16}}
\put(154,41){\footnotesize{3}}
\put(164,44){\line(1,0){8}}
\put(176,44){\circle*{8}}
\put(173,50){\scriptsize{-3}}
\put(180,44){\line(1,0){8}}
\put(192,44){\circle*{8}}
\put(189,50){\scriptsize{-4}}
\put(196,44){\line(1,0){8}}
\put(212,44){\circle{16}}
\put(210,41){\footnotesize{4}}
\put(220,44){\line(1,0){8}}
\put(232,44){\circle*{8}}
\put(229,50){\scriptsize{-2}}
\put(236,44){\line(1,0){8}}
\put(248,44){\circle*{8}}
\put(245,50){\scriptsize{-6}}
\put(252,44){\line(1,0){8}}
\put(268,44){\circle{16}}
\put(266,41){\footnotesize{5}}
\put(276,44){\line(1,0){8}}
\put(288,44){\circle*{8}}
\put(286,50){\scriptsize{-2}}
\put(292,44){\line(1,0){8}}
\put(304,44){\circle*{8}}
\put(301,50){\scriptsize{-3}}
\put(308,44){\line(1,0){8}}
\put(320,44){\circle*{8}}
\put(317,50){\scriptsize{-2}}
\put(324,44){\line(1,0){8}}
\put(336,44){\circle*{8}}
\put(333,50){\scriptsize{-2}}
\put(28,37){\line(0,-1){20}}
\put(28,9){\circle{16}}
\put(26,6){\footnotesize{6}}
\put(36,9){\line(1,0){8}}
\put(48,9){\circle*{8}}
\put(45,15){\scriptsize{-2}}
\put(52,9){\line(1,0){8}}
\put(64,9){\circle*{8}}
\put(61,15){\scriptsize{-2}}
\put(68,9){\line(1,0){8}}
\put(80,9){\circle*{8}}
\put(77,15){\scriptsize{-2}}
\put(84,9){\line(1,0){8}}
\put(96,9){\circle*{8}}
\put(93,15){\scriptsize{-2}}
\put(100,9){\line(1,0){8}}
\put(112,9){\circle*{8}}
\put(109,15){\scriptsize{-3}}
\put(116,9){\line(1,0){8}}
\put(132,9){\circle{16}}
\put(130,6){\footnotesize{7}}
\put(140,9){\line(1,0){8}}
\put(152,9){\circle*{8}}
\put(149,15){\scriptsize{-3}}
\put(156,9){\line(1,0){8}}
\put(168,9){\circle*{8}}
\put(165,15){\scriptsize{-4}}
\put(172,9){\line(1,0){8}}
\put(188,9){\circle{16}}
\put(186,6){\footnotesize{8}}
\put(196,9){\line(1,0){8}}
\put(208,9){\circle*{8}}
\put(205,15){\scriptsize{-2}}
\put(212,9){\line(1,0){8}}
\put(224,9){\circle*{8}}
\put(221,15){\scriptsize{-6}}
\put(228,9){\line(1,0){8}}
\put(244,9){\circle{16}}
\put(242,6){\footnotesize{9}}
\put(252,9){\line(1,0){8}}
\put(264,9){\circle*{8}}
\put(261,15){\scriptsize{-2}}
\put(268,9){\line(1,0){8}}
\put(280,9){\circle*{8}}
\put(277,15){\scriptsize{-3}}
\put(284,9){\line(1,0){8}}
\put(296,9){\circle*{8}}
\put(293,15){\scriptsize{-2}}
\put(300,9){\line(1,0){8}}
\put(312,9){\circle*{8}}
\put(309,15){\scriptsize{-2}}
\put(316,9){\line(1,0){8}}
\put(332,9){\circle{16}}
\put(327,6){\footnotesize{10}}
\put(332,1){\line(0,-1){8}}
\put(332,-11){\circle*{8}}
\put(314,-14){\scriptsize{-11}}
\put(332,-15){\line(0,-1){8}}
\put(332,-31){\circle{16}}
\put(327,-34){\footnotesize{11}}
\put(264,-31){\circle*{8}}
\put(261,-25){\scriptsize{-2}}
\put(268,-31){\line(1,0){8}}
\put(280,-31){\circle*{8}}
\put(277,-25){\scriptsize{-3}}
\put(284,-31){\line(1,0){8}}
\put(296,-31){\circle*{8}}
\put(293,-25){\scriptsize{-2}}
\put(300,-31){\line(1,0){8}}
\put(312,-31){\circle*{8}}
\put(309,-25){\scriptsize{-2}}
\put(316,-31){\line(1,0){8}}
\end{picture}
\end{center}

\vspace{0.2cm}
\begin{corollary} In the case $55(0)$ we have $I=11$, $1\le \rho(S)\le 11$ and
$\rank\Delta(S)+\rho(S)=12$.
\end{corollary}

\begin{theorem} In the case $56(0)$ the set $\TT$ must satisfy the following
condition: $1\in\TT_1$.
\end{theorem}

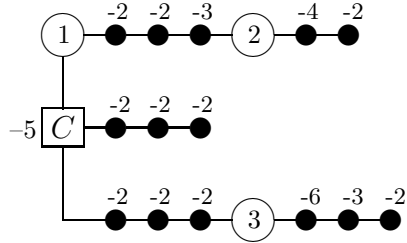
\begin{center}
\begin{picture}(156,90)(0,0)
\put(8,40){\footnotesize{--5}}
\put(20,40){\fbox{$C$}}
\put(28,51){\line(0,1){20}}
\put(28,79){\circle{16}}
\put(26,76){\footnotesize{1}}
\put(36,79){\line(1,0){8}}
\put(48,79){\circle*{8}}
\put(45,85){\scriptsize{-2}}
\put(52,79){\line(1,0){8}}
\put(64,79){\circle*{8}}
\put(61,85){\scriptsize{-2}}
\put(68,79){\line(1,0){8}}
\put(80,79){\circle*{8}}
\put(77,85){\scriptsize{-3}}
\put(84,79){\line(1,0){8}}
\put(100,79){\circle{16}}
\put(98,76){\footnotesize{2}}
\put(108,79){\line(1,0){8}}
\put(120,79){\circle*{8}}
\put(117,85){\scriptsize{-4}}
\put(124,79){\line(1,0){8}}
\put(136,79){\circle*{8}}
\put(134,85){\scriptsize{-2}}
\put(36,44){\line(1,0){8}}
\put(48,44){\circle*{8}}
\put(46,50){\scriptsize{-2}}
\put(52,44){\line(1,0){8}}
\put(64,44){\circle*{8}}
\put(61,50){\scriptsize{-2}}
\put(68,44){\line(1,0){8}}
\put(80,44){\circle*{8}}
\put(77,50){\scriptsize{-2}}
\put(28,37){\line(0,-1){28}}
\put(28,9){\line(1,0){16}}
\put(48,9){\circle*{8}}
\put(45,15){\scriptsize{-2}}
\put(52,9){\line(1,0){8}}
\put(64,9){\circle*{8}}
\put(61,15){\scriptsize{-2}}
\put(68,9){\line(1,0){8}}
\put(80,9){\circle*{8}}
\put(77,15){\scriptsize{-2}}
\put(84,9){\line(1,0){8}}
\put(100,9){\circle{16}}
\put(98,6){\footnotesize{3}}
\put(108,9){\line(1,0){8}}
\put(120,9){\circle*{8}}
\put(118,15){\scriptsize{-6}}
\put(124,9){\line(1,0){8}}
\put(136,9){\circle*{8}}
\put(134,15){\scriptsize{-3}}
\put(140,9){\line(1,0){8}}
\put(152,9){\circle*{8}}
\put(150,15){\scriptsize{-2}}
\end{picture}
\end{center}

\vspace{0.2cm}
\begin{corollary} In the case $56(0)$ we have $I=7$, $1\le \rho(S)\le 3$ and
$\rank\Delta(S)+\rho(S)=13$.
\end{corollary}

\end{document}